\documentclass{article}
\usepackage{amsmath,amssymb,amsthm}
\usepackage{atbegshi}
\AtBeginDocument{\AtBeginShipoutNext{\AtBeginShipoutDiscard}}
\usepackage{graphicx}
\usepackage{color}
\usepackage[english]{babel}
\usepackage[affil-it]{authblk}
\usepackage{bm}
\numberwithin{equation}{section}
\theoremstyle{plain}

\newtheorem{lemma}{Lemma}[section]

\begin{document}
\title{Stepwise Methods in Optimal Control Problems}

\author{Mehdi Afshar}
 \affil{ Department of Mathematical Sciences, Institute for Advanced Studies in basic Sciences,\\ GavaZang, Zanjan, IRAN}
 \thanks{\texttt{m\_afshar@iasbs.ac.ir}}
  
\author{Farshad Merrikhbayat}
\affil{Department of Electrical and Computer Engineering, Zanjan University,\\ P.O.Box 313, Zanjan, IRAN}
\thanks{\texttt{f.bayat@znu.ac.ir}}

\author{Mohammad Reza Razvan}
\affil{Department of Mathematical Sciences, Sharif University of Technology,\\
  P.O. Box 11155-9415, Tehran, IRAN}
\thanks{\texttt{razvan@sina.sharif.edu}}

\maketitle

\begin{abstract}
We introduce a new method, stepwise method for solving optimal control problems. Our first motivation for new approach emanate from limitations on continuous time control functions in PMP. Practically in most of the real world models, we are not able to change control value for every time such as in drug dose calculation or in resourse allocation problems. But it is practical to change control value in some time section that lead to stepwise control function. We study some examples via classical Pontryagin Maximum Principle(PMP) and via stepwise method.  The new method has some other advantages in comparison with PMP method in models with complicated cost function or systems. In real world applications, the new method has a high performance in implementation.
\end{abstract}

\section{Introduction}
Optimal control theory is an effective tool in real world modelling such as physical, biological, economical and other models. Diverse examples are studied in\cite{sethi} and \cite{lenhart} .  Optimal control theory is used in chemotherapy of cancer\cite{fisterpan}. Several papers are studied about epidemiology from optimal control theory viewpoint, for example \cite{DSDI}.
The most important classical method in optimal control theory is the remarkable result of Pontryagin Maximum Principle(PMP) that be used in various forms in applied problems. there are necessary conditions in PMP that lead to limitations in applications. These conditions can occurred for functions in state equations, cost function and control functions. Here we interested in finding a new optimal control method with more  proficiency in complicated cases. In this technique, control functions selected among the stepwise functions. Another notable point  in the new manner is the combination of the heuristic and classical methods. In the forthcoming sections, numerical forward-backward sweep method and stepwise method are applied to some problems such that one can obtain a clear vision about potency and power of the new method. Furthermore stepwise method works easily in the problems with more complicated cost functions in contrast to classical methods. In this paper the proficiency of stepwise method is shown in some models.

\section{ Introductory example and definitions of stepwise method with fixed step-size}
In this section, we describe the stepwise method through the simple example. Consider the problem
\begin{center}
$max\{J=\int_{0}^{2}(2x-3u-u^2)dt\}$
\end{center}
subject to $\dot{x}=x+u,\quad x(0)=5$ and the control constrain $u\in \Omega=[0,2].$
Solution: 
We can find the optimal control function $u(t)$ via PMP. The Hamiltonian is 
\begin{center}
$H=(2x-3u-u^2)+\lambda(x+u)=(2+\lambda)x-(u^2+3u-\lambda u).$
\end{center}
One can find the optimal control policy by differentiating $H$ with respect to $u$. Thus 
\begin{center}
$\frac{\partial H}{\partial u}=-2u-3+\lambda=0,$
\end{center}
so that the control function is $u(t)=\frac{\lambda(t)-3}{2}$ that $u(t)$ stays within the interval $\Omega=[0,2].$ We next drive the adjoint equation as 
\begin{equation*}
\begin{aligned}
&\dot{\lambda}=-\frac{\partial H}{\partial x}=-2-\lambda,\quad \lambda(2)=0\\&\dot{\lambda}+\lambda=-2, \quad \lambda(2)=0.
\end{aligned}
\end{equation*}
This equation can be solved and $\lambda(t)=2(e^{2-t}-1)$ when we impose the control constraint $\Omega=[0,2],$ the optimal control is obtained:
\begin{equation}\nonumber
u=\left\{\begin{array}{ccccc}
               2 & \mbox{if \quad $e^{2-t}-2.5>2$},   \\
               e^{2-t}-2.5 &\mbox{if \quad$0\leq e^{2-t}-2.5\leq 2$},  \\
               0 &\mbox{if \quad$e^{2-t}-2.5<0$}.
          \end{array}\right.
\end{equation}
\begin{figure}[htbp]
    \begin{center}
    \scalebox{0.35}{\includegraphics{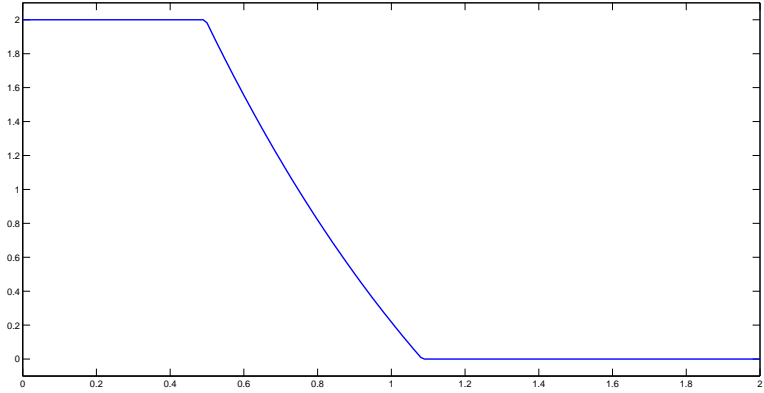}}
    \caption{Graph of optimal control for introductory example.}
    \label{fig:seth}
    \end{center}
\end{figure}
 The final cost $J$ is $68.93$. For solving this problem by new method, we changed the maximization problem to minimization through converting $J$ to $1/(1+J)$ and the final cost in minimization problem is $0.0143$.
When we use PMP, the control functions must satisfy the special conditions. In practical applications, we are not able to change control value every moment continuously. But we can change control values in some slices of time that lead to step function. In order to solve this problem by stepwise method, we seek control function among stepwise functions. For this purpose, we divide $[0,T]$ into equal parts. Suppose that the control function $u(t)$ has a constant value in each part. 

Let us enter the stepwise control function in control system
{
    \def\OldComma{,}
    \catcode`\,=13
    \def,{%
      \ifmmode%
        \OldComma\discretionary{}{}{}%
      \else%
        \OldComma%
      \fi%
    }%
$\dot{x}=f(x,u,t), x(0)=x_0$.
 }
 For $t\in [0,\frac{T}{3}]$,, we have $u(t)=\alpha$ and one can solve the ode $\dot{x}=f(x,\alpha,t),\quad x(0)=x_0$.  For $t\in [\frac{T}{3},\frac{2T}{3}]$, we solve the ode $\dot{x}=f(x,\beta,t)$ together with initial condition $ x(0)=x(\frac{T}{3})$ that is the terminal point of system in $t\in [0,\frac{T}{3}]$. Same procedure repeat for $\dot{x}=f(x,\gamma,t),\quad x(0)=x(\frac{2T}{3})$ when $t\in [\frac{2T}{3},T]$. We can compute cost for a typical stepwise control function and search for the optimal values for $(\alpha,\beta,\gamma)$. In this manner we convert the optimal control problem to optimization problem. For solving the optimization problems, we are able to use some analytical and heuristic and metaheuristic methods such as pattern search, simulated annealing, genetic algorithm, and other methods which these methods works easily in the problems with complicated systems and cost functions. There may be exist an important question here about the final cost. Is it possible that the difference between final cost in PMP method and in stepwise method exceeds from our expectation There is a simple lemma about stepwise functions and continuous functions can be helpful here. 
\begin{lemma}
For every continuous function $u(t)$, there is a sequence $\{u_n(t)\}$ of stepwise functions that $\lim_{n\to\infty} u_n(t)= u(t).$
\end{lemma}
Using this lemma, we can be confident that the new method do not generate useless solutions. Note that some step functions are not able to satisfy the PMP condition and do not belong to admissible controls. But it is likely, we can find step function in new method with lower cost than PMP solution. The final cost in the stepwise method equal to $0.014305429952056222$. We illustrate this method in figure below. Note that ,We have obtained  all the numerical results in 20-30 run times.
 \begin{figure}[htbp]
    \begin{center}
    \scalebox{0.35}{\includegraphics{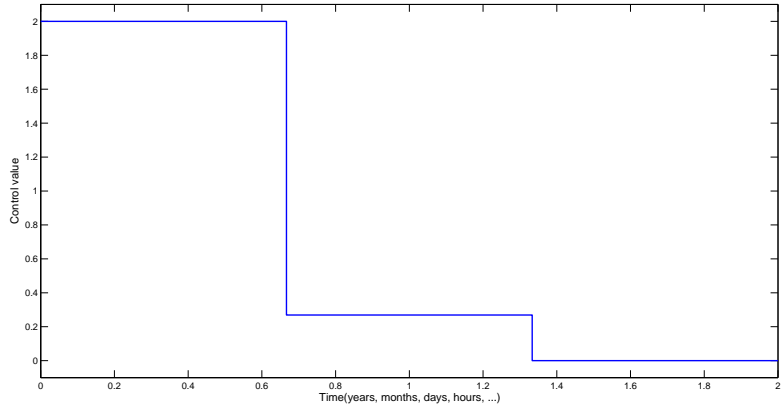}}
    \caption{Graph of optimal control for introductory example via stepwise method.}
    \label{fig:stair}
    \end{center}
\end{figure}
We continue our study on 3-step stepwise method via pattern search, simulated annealing and genetic algorithm. The results is illustrated in next figure and table.
\begin{center}
\begin{tabular}{ |l|l| }
    \hline
     Method & Final cost (J) \\ \hline
 Pattern search &  $0.014305429952056222$ \\
 Simulated annealing & $0.014306824417196181$\\
 Genetic algorithm & $0.014358890512889641$ \\
\hline
\end{tabular}
\end{center}
 \begin{figure}[htbp]
    \begin{center}
    \scalebox{0.35}{\includegraphics{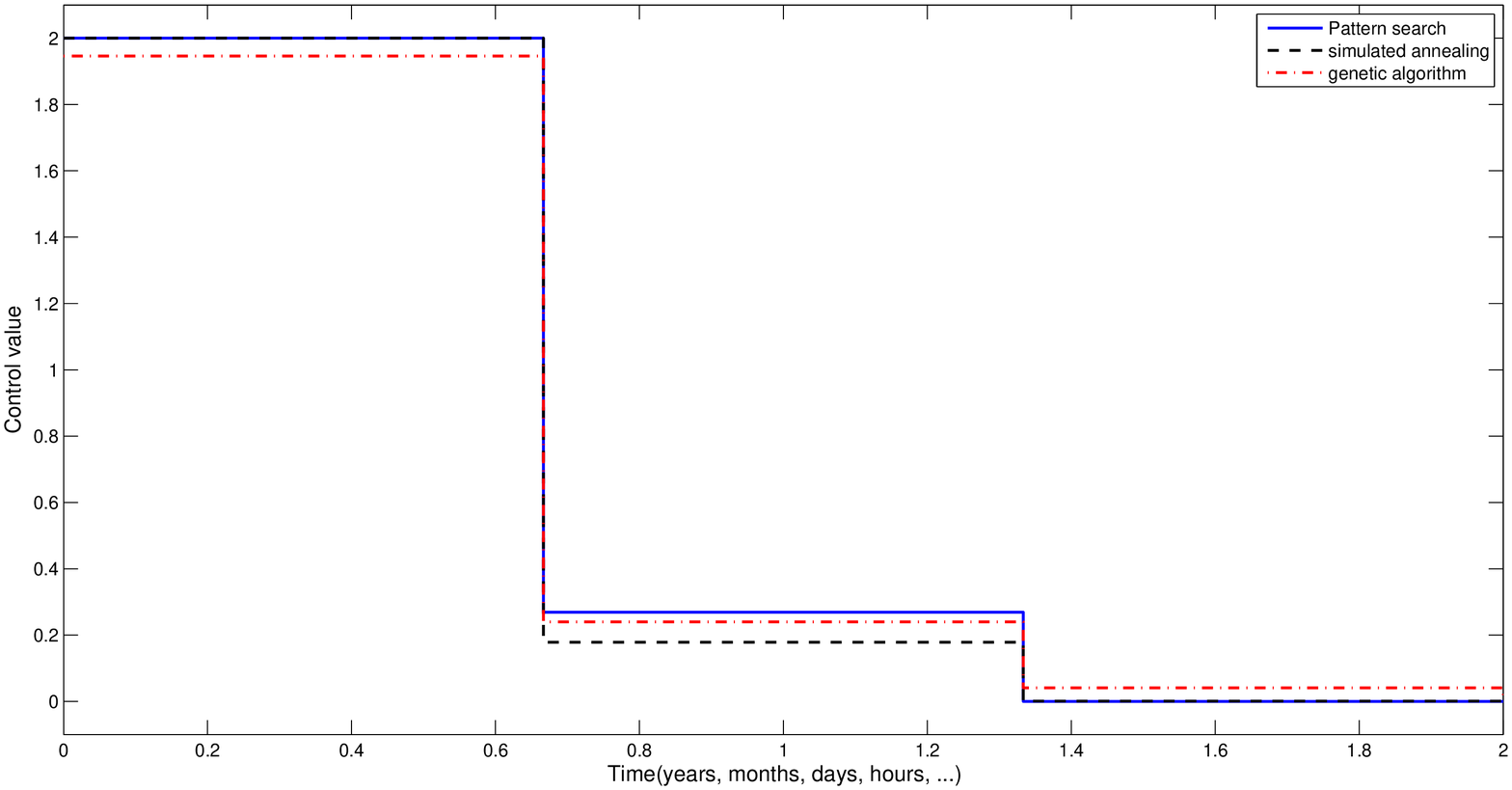}}
    \caption{Results for stepwise method via pattern search, simulated annealing and genetic algorithm.}
    \label{fig:stair}
    \end{center}
\end{figure}
 Trying to get better stepwise control function through additional steps seems natural. We applied 5-step function instead 3-step function.These results coincide with our expectations.
\begin{center}
\begin{tabular}{ |l|l| }
    \hline
     Method & Final cost (J) \\ \hline
 Pattern search &  $0.014283994290191705$ \\
 Simulated annealing & $0.01430531425509952$\\
 Genetic algorithm & $0.0143199869929726 $ \\
\hline
\end{tabular}
\end{center}

 \begin{figure}[htbp]
    \begin{center}
    \scalebox{0.35}{\includegraphics{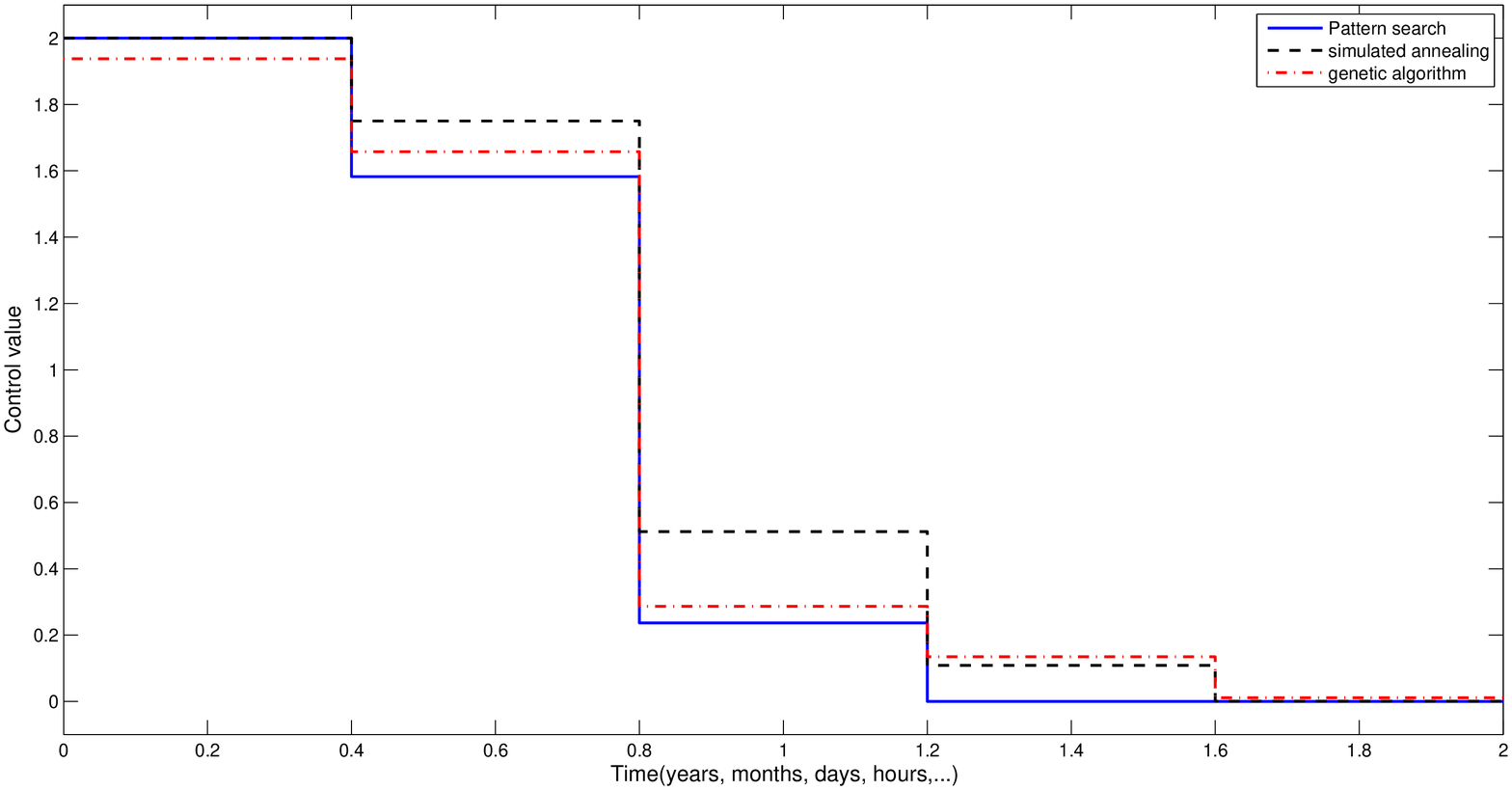}}
    \caption{Results for 5-step stepwise method via pattern search, simulated annealing and genetic algorithm.}
    \label{fig:stair}
    \end{center}
\end{figure}
\section{Stepwise method with variable step-size }
In the former section, we divided the interval into some equal parts. Here we let the optimization method decide about width of subintervals. Let's come back to introductory example. Here, instead of dividing $[0,2]$ to $[0,1/3]$,$[1/3,2/3]$ and $[2/3,2]$, we divide $[0,a]$,$[a,b]$ and $[b,2]$ and let the optimization method decide about $a$ and $b$. The next figures and tables  show the improvement in optimal policy and final cost and give complete information about subinterval and control value on the subinterval in variable stepwise method.

\begin{center}
\begin{tabular}{ |l|l|}
    \hline
     Method & Final cost (J)  \\ \hline
 Pattern search &  $0.012566299700335069$ \\
 Simulated annealing & $0.01334206413715155$\\
 Genetic algorithm & $0.012912036300546229 $\\
\hline
\end{tabular}
\end{center}
\small
\begin{center}
\begin{tabular}{ |l|l|l|}
    \hline
   Method& Subintervals  & control value \\ \hline
 Pattern search    &[0,0],[0,1],[1,2]&(0,2,0) \\
 Simulated annealing &[0,0.0036],[0.0036,0.9738],[0.9738, 2]&(1.6336,1.8345,0.5623)\\
 Genetic algorithm &[0,0.0034],[0.0034,0.9027],[0.9027,2]&(0.7718,1.9087,0.1524) \\
\hline
\end{tabular}
\end{center}
\normalsize
 \begin{figure}[htbp]
    \begin{center}
    \scalebox{0.3}{\includegraphics{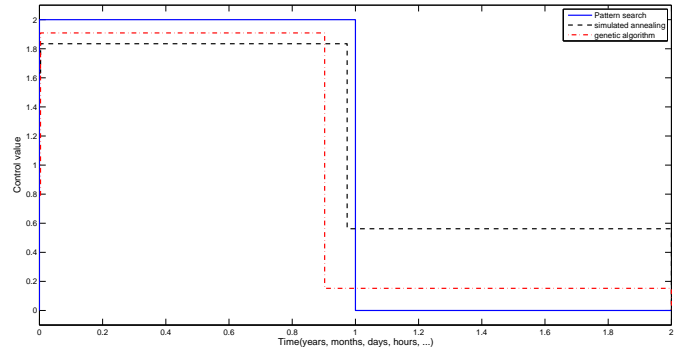}}
    \caption{Results for variable 3-step stepwise method via pattern search, simulated annealing and genetic algorithm.}
    \label{fig:stair}
    \end{center}
\end{figure} 
\section{Stepwise method in real world models}
As we mentioned before, there are limitations for admissible controls in PMP approach such as continuity with respect to time and others. Practically, we are not able to change the value of control function in every moment of time interval. Instead, one can change the control value at several time sections.Thus, it seems that the stepwise method is a reasonable way in real world applications.  For example, when you want to make a decision about resource allocation in epidemiological models, you can not alter your strategy in short time interval. Because changing the vaccination rate or prevention strategy may imposes heavy costs. There are same problems in optimal control model of treatment  of disease through the use of drugs. It seems that there are sufficient motivation to practice the stepwise method in real world processes. The next examples show the performance of stepwise method in contrast to classical PMP method.
\subsection{Example: Chemotherapy }
Optimal control methods are useful in optimal control model of chemotherapy. Renee fister et al in\cite{fisterpan} studied different cell-kill models of chemotherapy. They characterized optimal control strategy that minimizes the cancer mass and the cost of total amount of drug. We use stepwise method for some of their models. The problem is:
\begin{equation*}
\min_u \int_{0}^{T}a(N(t)-N_d)^2+bu^2(t) dt
\end{equation*} 
subject to,
\begin{equation*}
\begin{aligned}
& N'(t)=rNln(\frac{1}{N})-u(t)\delta N(t)\\& N(0)=N_0,\quad u(t)\geq 0.
\end{aligned}
\end{equation*} 
The following parameters appear in model:\\
$N(t)$: The normalized density of the tumor at time $t$\\
$r$: The growth rate of the tumor\\
$\delta$: The magnitude of the dose\\
$u(t)$: The time dependent pharmacokinetics of the drug \\
$N_d$: The desired tumor density.\\ Let us enter these value: $r=0.1$, $a=3$, $b=1$, $\delta=0.45$, $N_d=0$, $N_0=0.975$ and $T=20$.
The next figures show the optimal control strategy in PMP and stepwise method.  
 \begin{figure}[htbp]
    \begin{center}
    \scalebox{0.3}{\includegraphics{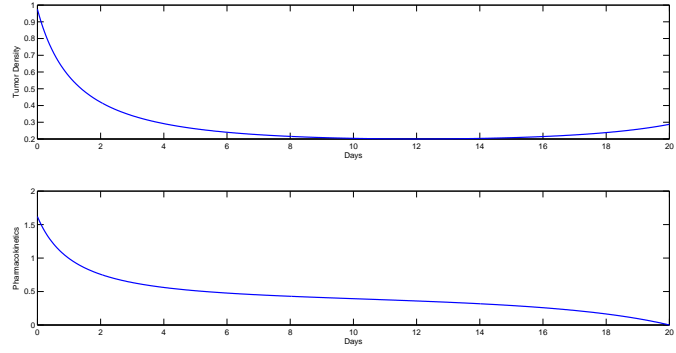}}
    \caption{Graph of optimal control for chemotherapy example via PMP method.}
    \label{fig:stair}
    \end{center}
\end{figure}
Now, we try to use stepwise method in this model and present the results below.
 \begin{figure}[htbp]
    \begin{center}
    \scalebox{0.3}{\includegraphics{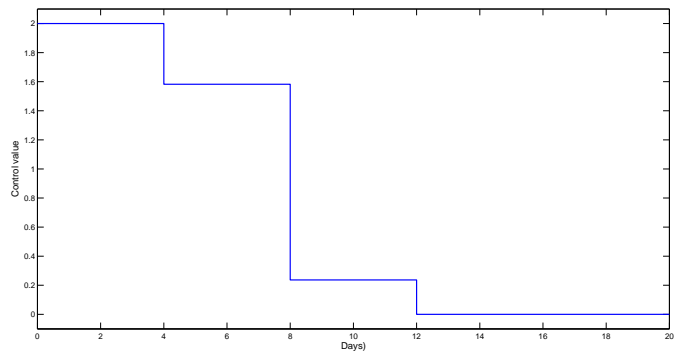}}
    \caption{Graph of optimal control for introductory example.}
    \label{fig:stair}
    \end{center}
\end{figure}
The final cost in PMP method equal to $10.7758$ and in 5-step stepwise method is $10.866632710096287$.
\vspace{-.3cm}
\subsection{Example: Differential susceptibility and differential infectivity model}
Based on \cite{DSDI}, we develop optimal control formulation of DSDI model with two groups of susceptible and two groups of infected individuals\cite{DSDIme}. Because of apparent diversity of examples, the
idea of dividing susceptible and infected population into two subgroups examined. For
example in plenty of diseases, disease processes is different in male or
female, children or adult, adicted or nonaddicted, and so on. The following parameters appear in our model:\\
$\mu$: natural death rate;\\
$\nu_i$: the rate at which infectives in $I_i$ are removed or become immune;\\
$\delta$: disease-induced mortality rates for the infectives;\\
$\lambda_i$: The rate of infection for susceptibles in group Si; \\
The infectivity rate $\lambda_i$ is given by $\lambda_i=r\alpha_i\sum_{j=1}^{2}\beta_jI_j$ in which $'\beta_i'$ is the transmission probability per contact and $'r'$ is the number of contacts of an individual per unit time.
 We suggested the following ODEs system (\ref{system2}) describing the model with controls.
\begin{equation}\label{system2}
\begin{cases}
\dot{S_1}&= \mu (p_1S^0- S_1)-\lambda_1 S_1(1-u_1)\\
\dot{S_2}&= \mu (p_2S^0- S_2)-\lambda_2 S_2(1-u_2)\\
\dot{I_1}&= q_{11}\lambda_1 S_1(1-u_1)+q_{21}\lambda_2 S_2(1-u_2)-(\mu+\nu_1+u_3)I_1\\
\dot{I_2}&= q_{12}\lambda_1 S_1(1-u_1)+q_{22}\lambda_2 S_2(1-u_2)-(\mu+\nu_2+u_4)I_2\\
\dot{R}&=(\nu_1+u_3)I_1+(\nu_2+u_4)I_2-(\mu+\delta)R
\end{cases}
\end{equation}
The control functions $u_1(t)$, $u_2(t)$, $u_3(t)$ and $u_4(t)$ have to be bounded on $[0,1]$ and Lebesgue integrable
functions. $u_1(t)$ and $u_2(t)$ measure the time dependent efforts on the preventive strategy of
susceptible individuals in $S$, to reduce the number of individuals that may be infectious.
The control functions $u_2(t)$ and $u_3(t)$ measures the time dependent efforts on the treatment of
infected individuals in $I_1$ and $I_2$ respectively. This control will have an impact on the output flow of people from the
The objective functional to be minimized is:\\
\begin{equation}\label{objective}
J(u_1,u_2,u_3,u_4)=\int^{T}_{0}{AI_1^2+BI_2^2+Cu_1^2+Du_2^2+Eu_3^2+Fu_4^2}dt
\end{equation}
Here, $A,B,C,D,E$ are adjustment parameters. We seek an optimal control triple $(u_1^*,u_2^*,u_3^*,u_4^*)$ such that
\begin{equation*}
J(u_1^*,u_2^*,u_3^*,u_4^*)=\min{\{J(u_1,u_2,u_3,u_4) | (u_1,u_2,u_3,u_4)\in U}\}
\end{equation*}
where
\small
$U=\{J(u_1,u_2,u_3,u_4) | u_i \mbox{ measerable} , 0\leq u_i \leq 1 , t\in [0,T], i=1,2,3,4\}$
\normalsize
 is the control set. Let us enter the following values in model system.
 \small
\begin{center}
\begin{tabular}{ |l|l|l| }
\hline
\multicolumn{2}{|c|}{Parameters and values} \\
\hline
 $S^0=1$ &  $\delta=0$ \\
 $\mu=.012$ & $S_1(0)=0.47$\\
 $T=1000$&$S_2(0)=0.47$\\
 $p_1=0.5$ & $I_1(0)=0.02$\\
 $p_2=0.5$ & $I_2(0)=.04$\\
 $\alpha_1=0.05$ & $R(0)=0$ \\
 $\alpha_2=0.2$&$\beta_1=0.2$ \\
 $\nu_1=0.15$ & $\beta_2=0.06$ \\
 $\nu_2=0.6$ &$r=25 $ \\
 $q_{11}=0.9$&$q_{12}=0.1$ \\
 $q_{21}=.1$& $q_{22}=.9$  \\
 $A=3$ & $B=3$  \\
 $C=0.002$ & $D=0.002$ \\
 $E=0.002$ & $F=0.002$ \\
\hline
\end{tabular}
\end{center}
\normalsize
Below, we can depict optimal control policy by PMP method and stepwise method(with pattern search for optimization problem).The final cost in PMP method equal to $0.1059$ and in 5-step stepwise method is $0.11107136532373643 $.
 \begin{figure}[htbp]
    \begin{center}
    \scalebox{0.31}{\includegraphics{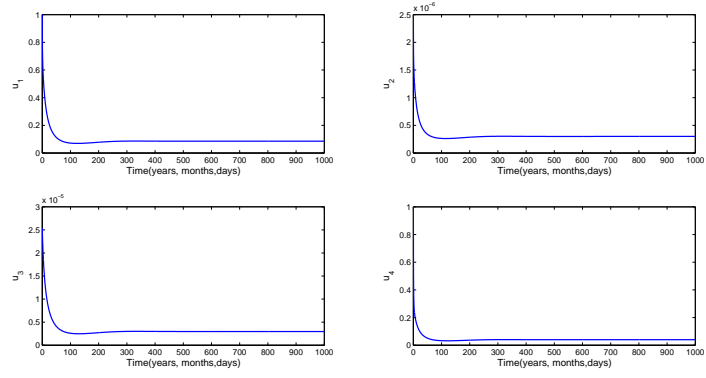}}
    \caption{Graph of optimal control for DSDI model via PMP.}
    \label{fig:stair}
    \end{center}
\end{figure}
 \begin{figure}[htbp]
    \begin{center}
    \scalebox{0.35}{\includegraphics{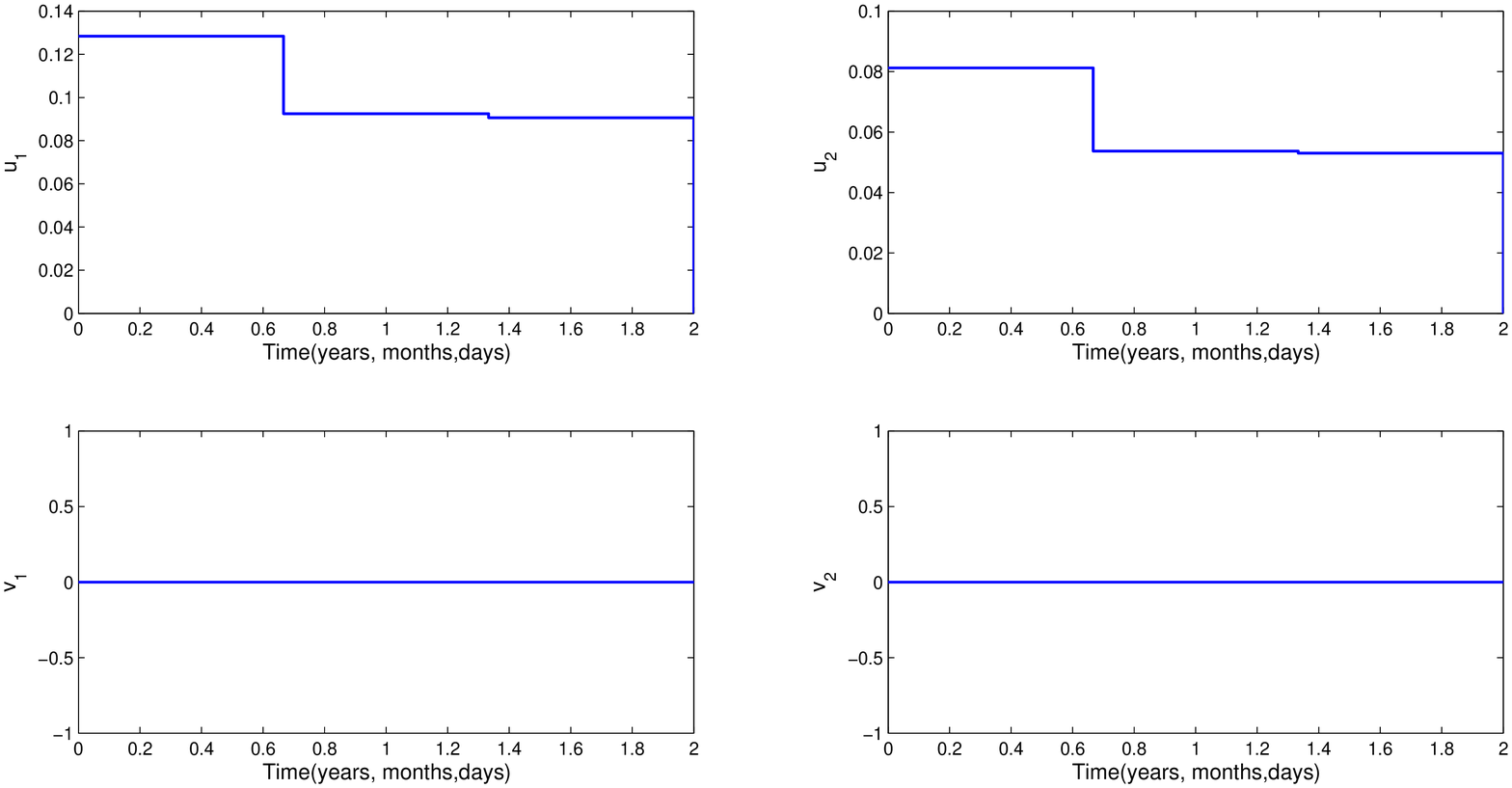}}
    \caption{Graph of optimal control DSDI model via 3-step stepwise method.}
    \label{fig:stair}
    \end{center}
\end{figure}

\section{conclusion}
we introduce the stepwise method for optimal control problems. This method could be replaced with PMP classic method in real world problems. Using this new method in various cases of applied models seems reasonable .
\section{Acknowledgement}
It is a pleasure to acknowledge the helpful suggestions made by Dr rooin during the preparation of this paper.

\end{document}